\documentclass[11pt]{amsart}
 \newtheorem{thm}{Theorem}[section]
 \newtheorem{cor}[thm]{Corollary}
 
 \newtheorem{prop}[thm]{Proposition}
 \theoremstyle{definition}
 \newtheorem{defn}[thm]{Definition}
 \theoremstyle{remark}
 
 \numberwithin{equation}{section}
 \newcommand{\norm}[1]{\left\Vert#1\right\Vert}
 
 \newcommand{\C}{\mathbb{C}}

\begin{document}

\title[Every transcendental
 operator has a non-trivial invariant subspace]
 {Every transcendental
 operator has a non-trivial invariant subspace}

\author{Yun-Su Kim  }

%\address{Department of Mathematics, Indiana University, Bloomington,
%Indiana, U.S.A. }

\address{Department of Mathematics, University of Toledo, Toledo,
Ohio, U.S.A. } \email{Yun-Su.Kim@utoledo.edu}


\dedicatory{}

%\commby{Daniel J. Rudolph}

%%% ----------------------------------------------------------------------

\newpage
%%% ----------------------------------------------------------------------

\begin{abstract}In this paper, to solve the invariant subspace problem, contraction operators are classified into
three classes ; (Case 1) completely non-unitary contractions with
a non-trivial algebraic element, (Case 2) completely non-unitary
contractions without a non-trivial algebraic element, or (Case 3)
contractions which are not completely non-unitary.

We know that every operator of (Case 3) has a non-trivial
invariant subspace. In this paper, we answer to the invariant
subspace problem for the operators of (Case 2). Since (Case 1) is
simpler than (Case 2), we leave as a question.

\end{abstract}

%%% ----------------------------------------------------------------------
\maketitle
%%% ----------------------------------------------------------------------

\section*{Introduction}
An important open problem in operator theory is the invariant
subspace problem. The invariant subspace problem is the question
whether the following statement is true or not: \vskip0.2cm

 Every bounded linear operator $T$ on a
separable Hilbert space $H$ of dimension $\geq{2}$ over $\C$ has a
non-trivial invariant subspace. \vskip0.2cm

Since the problem is solved for all finite dimensional complex
vector spaces of dimension at least 2, in this note, $H$ denotes a
separable Hilbert space whose dimension is infinite. It is enough
to think for a contraction $T$, i.e., $\norm{T}\leq{1}$ on $H$.
Thus, in this note, $T$ denotes a contraction.

If $T$ is a contraction, then\vskip0.1cm

(Case 1) $T$ is a completely non-unitary contraction with a
\emph{non-trivial} \emph{algebraic element}, or\vskip.1cm

(Case 2) $T$ is a \emph{transcendental operator}; that is, $T$ is
a completely non-unitary contraction without a non-trivial
algebraic element, or\vskip.1cm

(Case 3) $T$ is not completely non-unitary.\vskip0.2cm

In this note, we discuss the invariant subspace problem for
operators of (Case 2).
 By using fundamental properties (Proposition \ref{2}, Proposition \ref{1}, and Corollary \ref{6}) of transcendental
 operators, we answer to the invariant subspace problem for the operators of (Case 2) in Theorem
 \ref{7}; \vskip0.1in

  \emph{Every transcendental
 operator defined on a separable Hilbert space $H$ has a non-trivial invariant subspace.}
\vskip0.1in Thus, we answered to the invariant subspace problem
for the (Case 2) in Theorem \ref{7} and, clearly, we know that
every operator of (Case 3) has a non-trivial invariant subspace.
Thus, to answer to the invariant subspace problem, it suffices to
answer for (Case 1).

We do not consider project (Case 1) in this note, and leave as a
question;

\vskip0.1cm

\textbf{Question}. Let $T\in{L(H)}$ be a completely non-unitary
contraction such that $T$ has a non-trivial algebraic element.
Then, does the operator $T$ have a non-trivial invariant subspace?

\vskip0.1cm The author would like to appreciate the advice of
Professor Ronald G. Douglas.

%\newpage
%%% ----------------------------------------------------------------------

%%% ----------------------------------------------------------------------
%%% ----------------------------------------------------------------------
%\centerline{\textbf{CHAPTER 1}}
%\vskip0.5cm\centerline{\textbf{Shift Operators}}

%\vskip0.5cm
%%% ----------------------------------------------------------------------
%\goodbreak
\section{Preliminaries and Notation}\label{13}

In this note, $\C$, $\overline{M}$ and $L(H)$ denote the set of
complex numbers, the (norm) closure of a set $M$, and the set of
bounded linear operators from $H$ to $H$ where $H$ is a separable
Hilbert space whose dimension is not finite, respectively.

For a set $A=\{a_{i}:i\in{I}\}\subset{H}$, $\bigvee{A}$ denotes
the closed subspace of $H$ generated by $\{a_{i}:i\in{I}\}$.

If $T\in{L(H)}$ and $M$ is an invariant subspace for $T$, then
$T|M$ is used to denote the restriction of $T$ to $M$, and
$\sigma(T)$ denotes the spectrum of $T$.

\subsection{A Functional Calculus.}\label{11}

 Let $H^{\infty}$
be the Banach space of all (complex-valued) bounded analytic
functions on the open unit disk $\textbf{D}$ with supremum norm
\cite{S2}. A contraction $T$ in $L(H)$ is said to be
\emph{completely non-unitary } provided its restriction to any
non-zero reducing subspace is never unitary.

%if there is no invariant subspace $K$ for $T$ such that $T|K$ is a
%unitary operator.

B. Sz.-Nagy and C. Foias introduced an important functional
calculus for completely non-unitary contractions.
\begin{prop}\label{12}Let $T\in{L(H)}$ be a completely non-unitary
contraction. Then there is a unique algebra representation
$\Phi_{T}$ from $H^{\infty}$ into $L(H)$ such that :\vskip0.2cm

(i) $\Phi_{T}(1)=I_{H}$, where $I_{H}\in{L(H)}$ is the identity
operator;

(ii) $\Phi_{T}(g)=T$, if $g(z)=z$ for all $z\in\textbf{D}$;

(iii) $\Phi_{T}$ is continuous when $H^{\infty}$ and $L(H)$ are
given the weak$^\ast$-

\quad\quad topology.

(iv) $\Phi_{T}$ is contractive, i.e.
$\norm{\Phi_{T}(u)}\leq\norm{u}$ for all
$u\in{H^{\infty}}$.\end{prop}

We simply denote by $u(T)$ the operator $\Phi_{T}(u)$.

B. Sz.- Nagy and C. Foias \cite{S2} defined the \emph{class
$C_{0}$} relative to the open unit disk \textbf{D} consisting of
completely non-unitary contractions $T$ on $H$ such that the
kernel of $\Phi_{T}$ is not trivial.  If $T\in{L(H)}$ is an
operator of class $C_{0}$, then \begin{center}$\ker
\Phi_{T}=\{u\in{H^{\infty}}:u(T)=0\}$\end{center} is a
weak$^{*}$-closed ideal of $H^{\infty}$, and hence there is an
inner function generating ker $\Phi_{T}$. The \emph{minimal
function} $m_{T}$ of an operator $T$ of class $C_{0}$ is the
generator of ker $\Phi_{T}$; that is, ker
$\Phi_{T}=m_{T}H^{\infty}$. Also, $m_{T}$ is uniquely determined
up to a constant scalar factor of absolute value one \cite{B1}.
%The theory
%of class $C_{0}$ relative to the open unit disk has been developed
%by B.Sz.- Nagy, C. Foias (\cite{S2}) and H. Bercovici (\cite{B1}).

\subsection{Algebraic Elements}
In this section, we provide the notion of \emph{algebraic
elements} for a completely non-unitary contraction $T$ in $L(H)$.

\begin{defn}\cite{K}
Let $T\in{L(H)}$ be a completely non-unitary contraction. An
element $h$ of $H$ is said to be \emph{algebraic with respect to}
$T$ provided that $\theta(T)h=0$ for some
$\theta\in{H^{\infty}}\setminus\{0\}$. If $h\neq{0}$, then $h$ is
said to be a \emph{non-trivial algebraic element} with respect to
$T$.

If $h$ is not algebraic with respect to $T$, then $h$ is said to
be \emph{transcendental with respect to} $T$.

\end{defn}

  % If $B$ is a closed subspace of $H$
%generated by $\{b_{i}\in{H}:i=1,2,3,\cdot\cdot\cdot\}$, then $B$
%will be denoted by $\bigvee_{n=1}^{\infty}b_{i}.$

\section{The Main Results}
%Let $T\in{L(H)}$ be a completely non-unitary contraction. To solve
%the invariant subspace problem, first, we prove that
%In this section, $H$ denotes a separable Hilbert space whose
%dimension is greater than or equal to two.
If $T$ is a contraction, then\vskip0.2cm

(Case 1) $T$ is a completely non-unitary contraction with a
non-trivial algebraic element, or

(Case 2) $T$ is a completely non-unitary contraction without a
non-trivial algebraic element; that is, every non-zero element in
$H$ is \emph{transcendental with respect to }$T$, or

(Case 3) $T$ is not completely non-unitary.\vskip0.2cm

It is clear for (Case 3). To answer to the invariant subspace
problem for (Case 2), we provide the following definition;

\begin{defn}
If $T$ is a completely non-unitary contraction without a
non-trivial algebraic element; that is, every non-zero element in
$H$ is transcendental with respect to $T$, then $T$ is said to be
a \emph{transcendental operator}.

\end{defn}

\begin{prop}\label{2} If $T:H\rightarrow{H}$ is a transcendental operator, then, for any
$\theta\in{H^{\infty}\setminus\{0\}}$, $\theta(T)$ is one-to-one.
\end{prop}
\begin{proof}
Suppose that $\theta(T)$ is not one-to-one for a function
$\theta\in{H^{\infty}\setminus\{0\}}$. Then, there is a non-zero
element $h$ in $H$ such that
\[\theta(T)h=0;\]
that is, $h$ is a non-trivial algebraic element with respect to
$T$. This, however, is a contradiction, since $T$ is a
transcendental operator. Thus, $\theta(T)$ is one-to-one for any
$\theta\in{H^{\infty}\setminus\{0\}}$.

\end{proof}

Recall that an arbitrary subset $M$ of $H$ is said to be
\emph{linearly independent} if every nonempty finite subset of $M$
is linearly independent.

\begin{prop}\label{1}
If $T:H\rightarrow{H}$ is a transcendental operator, then, for any
non-zero element $h$ in $H$,
$M=\{T^{n}h:n=0,1,2,\cdot\cdot\cdot\}$ is linearly independent.
\end{prop}
\begin{proof}
Let $h\in{H\setminus\{0\}}$ be given. Suppose that
$M=\{T^{n}h:n=0,1,2,\cdot\cdot\cdot\}$ is not linearly
independent. Then, there is a polynomial
$p\in{H^{\infty}}\setminus\{0\}$ such that $p(T)h=0$. Thus $h$ is
a non-trivial algebraic element with respect to $T$. This,
however, is a contradiction, since $T$ is a transcendental
operator. Therefore, $M$ is linearly independent.

\end{proof}

\begin{cor}\label{6}
Under the same assumption as Proposition \ref{1}, for a given
function $\theta\in{H^{\infty}\setminus\{0\}}$,
$M^{\prime}=\{\theta(T)^{n}h:n=1,2,\cdot\cdot\cdot\}$ is linearly
independent.

\end{cor}
\begin{proof}
In the same way as Proposition \ref{1}, it is proven.
\end{proof}

By a \emph{densely defined operator} $K$ in $H$, we mean a linear
mapping $K$ on the domain $\mathcal{D}(K)$ (which is a subspace of
$H$ and dense in $H$) of $K$ into $H$ \cite{16}. Recall that
\begin{equation}\label{24}
\mathcal{D}(SK)=\{x\in{\mathcal{D}(K)}:Kx\in{\mathcal{D}(S)}\},
\end{equation}
where $S$ and $K$ are unbounded operators \cite{16}.

Finally, the time has come to answer to the invariant subspace
problem for operators of (Case 2).

\begin{thm}\label{7}
If $T$ is  a transcendental operator in $L(H)$ and
$\lambda\in{\sigma(T)}$, then \vskip0.2cm

(i) $S=T-\lambda{I_{H}}$ has a non-trivial invariant subspace,
where $I_{H}$ is the identity operator on $H$, and \vskip0.1cm

(ii) $T$ also has a non-trivial invariant subspace.

\end{thm}

\begin{proof}

(i) Since $T$ is transcendental, $S\neq{0}$ (Note that if $S=0$,
then $p(T)=0$ where $p(z)=z-\lambda$), and by Proposition \ref{2},
$S(=p(T))$ is one-to-one. Since $S$ is not invertible, $S$ is not
onto. Let $h$ be a non-zero element in $H$ such that $h$ does not
belong to the range of $S$ and
\begin{equation}\label{3}M=\bigvee{\{S^{n}h:n=0,1,2,\cdot\cdot\cdot\}}.\end{equation}
If $M\neq{H}$, then $M$ is a non-trivial invariant subspace for
$S$, and so we assume that $M=H$. If
\begin{equation}\label{19}M^{\prime}=\bigvee{\{S^{n}h:n=1,2,\cdot\cdot\cdot\}},\end{equation}
then we will show that $h\notin{M^{\prime}}$. Since $S$ is
one-to-one and $h\neq{0}$,
\begin{equation}\label{70}M^{\prime}\neq\{0\}.\end{equation}

Suppose that $h\in{M^{\prime}}$.

By Corollary \ref{6}, we conclude that
$\{S^{n}h:n=1,2,\cdot\cdot\cdot\}$ is linearly independent, and
so, by Gram-Schmidt process, we have an orthonormal basis $B$ of
$M^{\prime}$ such that
\begin{equation}\label{21}B=\{P_{i}(S)h:i=1,2,\cdot\cdot\cdot\},\end{equation}
where $P_{i}(i=1,2,\cdot\cdot\cdot)$ is a polynomial satisfying
$\bigvee\{S^{n}h:n=1,2,\cdot\cdot,m\}=\bigvee\{P_{i}(S)h:i=1,2,\cdot\cdot,m\}$
for any $m\in\{1,2,\cdot\cdot\cdot\}$.

Then,
\begin{equation}\label{4}h=\sum_{i=1}^{\infty}a_{i}P_{i}(S)h\end{equation} where
$a_{i}\in{\C}$.

It follows that
\begin{equation}\label{5}\sum_{i=0}^{\infty}a_{i}P_{i}(S)h=0\end{equation}
where $a_{0}=-1$, and $P_{0}(S)=I_{H}$.

For any $k\in\{0,1,2,3,\cdot\cdot\cdot\}$, by equation (\ref{5}),
\begin{equation}\label{12}
\sum_{i=0}^{\infty}a_{i}P_{i}(S)(S^{k}h)=\lim_{m\rightarrow\infty}\sum_{i=0}^{m}a_{i}P_{i}(S)(S^{k}h)=
  \lim_{m\rightarrow\infty}S^{k}(\sum_{i=0}^{m}a_{i}P_{i}(S)h)=0
\end{equation}

By the same way as above, since
$\{S^{n}h:n=0,1,2,\cdot\cdot\cdot\}$ is linearly independent, by
Gram-Schmidt process, we have an orthonormal basis $B^{\prime}$ of
$H$ such that
\begin{equation}\label{21}B^{\prime}=\{e_{i}=f_{i}(S)h:i=0,1,2,\cdot\cdot\cdot\},\end{equation}
where $f_{i}(i=0,1,2,\cdot\cdot\cdot)$ is a polynomial satisfying
$\bigvee\{S^{n}h:n=0,1,2,\cdot\cdot,m\}=\bigvee\{f_{i}(S)h:i=0,1,2,\cdot\cdot,m\}$
for any $m\in\{0,1,2,\cdot\cdot\cdot\}$. Note that
$e_{0}=\frac{h}{\norm{h}}$.

Clearly, $\sum_{i=0}^{\infty}a_{i}P_{i}(S)$ is linear on
$\{c_{n}S^{n}h:c_{n}\in\C\texttt{ and }n=0,1,2,\cdot\cdot\cdot\}$,
and by equation (\ref{12}), $\sum_{i=0}^{\infty}a_{i}P_{i}(S)$ is
a densely defined operator in $H$. It is not assumed that
$\sum_{i=0}^{\infty}a_{i}P_{i}(S)$ is bounded or continuous.

By equations (\ref{12}) and (\ref{21}), we have that, for any
$e_{i}\in{B^{\prime}}$,
\begin{equation}\label{18}(\sum_{i=0}^{\infty}a_{i}P_{i}(S))e_{i}=0.\end{equation}
Thus, if
\begin{equation}\label{23}\mathcal{A}=\{\sum_{i=0}^{m}c_{i}e_{i}:c_{i}\in{\C},
\texttt{ and  }m=0,1,2,\cdot\cdot\cdot\},\end{equation} then, by
equation (\ref{18}),
\begin{equation}\label{22}(\sum_{i=0}^{\infty}a_{i}P_{i}(S))x\equiv{0}\end{equation}
for any
$x\in\mathcal{A}\subset\mathcal{D}(\sum_{i=0}^{\infty}a_{i}P_{i}(S))$.

 Let
\[K=\sum_{i=1}^{\infty}a_{i}g_{i}(S),\]
where $g_{i}$ is a polynomial such that
\begin{equation}\label{25}P_{i}(z)=zg_{i}(z)\end{equation} for
$i\in\{1,2,3,\cdot\cdot\cdot\}$. Note that by the definition of
$P_{i}$, we can easily find the polynomial $g_i$ satisfying
equation (\ref{25}) for $i\in\{1,2,3,\cdot\cdot\cdot\}$.

%Since $a_{0}=-1$, by equation
%(\ref{18}),
%\[K(h)=KS(h_{1})=\sum_{i=1}^{\infty}a_{i}S^{i}=h_{1},\texttt{ and }
%K(S^{n}h)=S^{n-1}h(n=1,2,3,\cdot\cdot\cdot).\]
Since $\{S^{n}h:n=1,2,3,\cdot\cdot\cdot\}\subset{\mathcal{D}(K)}$
( by equation (\ref{4})) and $h\in{M^{\prime}}$ (defined in
(\ref{19})) by assumption, $K$ is a densely defined operator in
$H$. It is not assumed that $K$ is bounded or continuous.

By equation (\ref{22}), since $a_{0}=-1$,
\begin{equation}\label{16}(SK)(x)=S(\sum_{i=1}^{\infty}a_{i}g_{i}(S))(x)
=(\sum_{i=1}^{\infty}a_{i}P_{i}(S))(x)=x,\end{equation} for any
$x\in\mathcal{A}\subset\mathcal{D}(K)$ (Note that
$\mathcal{A}\subset\mathcal{D}(SK)\subset{\mathcal{D}(K)}$ by
equation (\ref{24}) \cite{16}).

Since $e_{0}=\frac{h}{\norm{h}}$ and $e_{0}\in{\mathcal{A}}$, by
equation (\ref{16}), we have that $SK(e_{0})=\frac{h}{\norm{h}}$.
Thus, $h$ belongs to the range of $S$, but it is a contradiction.

Therefore,
\begin{equation}\label{71}h\notin{M^{\prime}}\end{equation}

Thus, by (\ref{70}) and (\ref{71}), we conclude that $M^{\prime}$
is a non-trivial invariant subspace for $S$. \vskip.2in

(ii) In $(i)$,
$T(M^{\prime})=(S+\lambda{I_{H}})M^{\prime}\subset{M^{\prime}}$

Therefore, $M^{\prime}$ is also a non-trivial invariant subspace
for $T$.

\end{proof}

%\begin{cor}

%In Theorem \ref{7}, we answered to the invariant subspace problem
%for the (Case 2). Finally, it is the time to answer to the
%invariant subspace problem;
%\begin{cor}
% Every bounded linear operator $T$ on a
%separable Hilbert space $H$ of dimension $\geq{2}$ over $\C$ has a
%non-trivial invariant subspace.

%\end{cor}

%If $T$ is a contraction, then\vskip0.2cm

%(Case 1) $T$ is a completely non-unitary contraction with a
%non-trivial algebraic element, or

%(Case 2) $T$ is a completely non-unitary contraction without a
%non-trivial algebraic element, or

%(Case 3) $T$ is not completely non-unitary.\vskip0.2cm

------------------------------------------------------------------------

\bibliographystyle{amsplain}
\bibliography{xbib}
\end{document}